\documentclass[11pt,a4paper]{article}
\usepackage[utf8]{inputenc}
\usepackage{amsmath, amsfonts, amssymb, amsthm}
\usepackage{geometry}
\usepackage{authblk}
\usepackage{tikz}
\newtheorem{theorem}{Theorem}[section]
\newtheorem{lemma}[theorem]{Lemma}
\newtheorem{proposition}[theorem]{Proposition}
\newtheorem{corollary}[theorem]{Corollary}

\newtheorem{definition}[theorem]{Definition}

\title{Non-commuting graph of AC-groups: as matroids}

\author[1]{Azizollah Azad\thanks{a-azad@araku.ac.ir}}
\author[2]{Nasim Karimi\thanks{nasim@ime.uerj.br}}
\author[1]{Sakineh Rahbariyan\thanks{s-rahbariyan@araku.ac.ir}}
\affil[1]{Department of Mathematics, Faculty of Science, Arak University, Arak, Iran}
\affil[2]{Instituto de Matem\'atica e Estat\'istica, Universidade do Estado de Rio de Janeiro, Brazil}

\date{\today}

\begin{document}

\maketitle

\begin{abstract}
Let $G$ be a finite non-abelian group with center $Z(G)$, and let $\Gamma_{G}$ denote its non-commuting graph, whose vertex set is $G \setminus Z(G)$ with edges connecting pairs of non-commuting elements \cite{AbdollahiAkbari}. In this paper, we characterize the class of AC-groups (groups with abelian centralizers) through the lens of matroid theory. By identifying $\Gamma_{G}$ with a 1-dimensional trim simplicial complex, we prove that $G$ is an AC-group if and only if the collection of its faces satisfies the matroid exchange axioms, rendering $\Gamma_G$ a matroid. This combinatorial framework bridges abstract group properties with geometric independence, providing a unified approach to studying commutativity in groups. Utilizing this structural duality, we derive an elegant formula for the total number of matroid bases and the clique number $\omega(\Gamma_{G})$ of AC-groups. Finally, we apply these results to compute the clique number for a specific family of CC-groups whose central quotient is a dihedral group, as well as certain families of finite $p$-groups, offering simplified, alternative proofs for several established results in the literature \cite{AzadRahbarian,FouladiOrfiAzad}.
\end{abstract}


\section{Introduction}
The study of graphs associated with algebraic structures has been a fruitful area of research for decades, tracing back to the foundational work on group-related graphs by Cayley \cite{Cayley} and later applications to finite groups by Bertram \cite{Bertram}. 

One of the most prominent examples explored in recent literature is the non-commuting graph $\Gamma_G$, whose vertex set consists of the non-central elements of a group, with edges connecting pairs of elements that do not commute \cite{AbdollahiAkbari}. The broad paradigm of studying graphs defined on algebraic structures has been extensively developed by Cameron in a series of highly influential surveys establishing a hierarchy of graphs on finite groups-- such as the power graph, enhanced power graph, deep commuting graph, and generating graph \cite{CameronSurvey, CameronSurvey2}. Within this rapidly expanding field, recent research has intensively investigated the combinatorial, algebraic, and spectral invariants of commuting and non-commuting graph variants for distinct algebraic structures, including structural coloring under automorphisms \cite{IlbeygiGanjali}, spectral energy characterizations of dihedral and dicyclic groups \cite{Assiry, Romdhini}, and extensions to Lie algebras \cite{Moghaddam}. The systematic study of these non-commuting graphs was significantly accelerated following Paul Erd\H{o}s' 1975 inquiry regarding the size of maximal subsets of pairwise non-commuting elements in a finite group \cite{Neumann}.

Matroid theory, initiated by Whitney \cite{Whitney} and Nakasawa \cite{NishimuraKuroda}, provides a general framework for independence that generalizes linear independence in vector spaces. While matroids are typically defined as abstract structures on a set $V$ satisfying the exchange property, they possess a specific characterization in the context of 1-dimensional simplicial complexes. Specifically, a graph $\Gamma$ functions as a matroid if and only if its complement $\overline{\Gamma}$ is a disjoint union of complete graphs. 
A group $G$ is called an AC-group if the centralizer $C_G(x)$ of every non-central element $x$ is abelian. In this paper, we bridge these two fields by proving a fundamental equivalence: a finite group $G$ is an AC-group if and only if its non-commuting graph $\Gamma_G$ is a matroid. Establishing this equivalence relation between AC-groups and matroids not only deepens our understanding of group structures but also provides a systematic algebraic machinery to generate matroids with highly interesting combinatorial and topological properties. In particular, this connection yields simple rank-2 matroids whose independence complexes exhibit well-defined face distributions, opening new avenues for exploring the interplay between finite groups and geometric lattices.

\section{Preliminaries}
In this section, we compile the necessary definitions, notation, and foundational results from group theory, matroid theory, and simplicial complexes that form the bedrock of our framework. For a deeper exposition of finite groups and non-commuting graphs, we refer to \cite{AbdollahiAkbari}. For standard terminology regarding matroids, we follow Oxley \cite{Oxley}, for graph theory we refer to Bondy and Murty \cite{BondyMurty}, and for simplicial complexes and geometric combinatorics, we refer to Bj\"orner \cite{Bjorner}.

\subsection{Group Theoretic Foundations and Non-commuting Graphs}
Throughout this paper, $G$ denotes a finite non-abelian group and $Z(G)$ denotes its center. 
\begin{definition}
For any non-central element $x \in G \setminus Z(G)$, the \textbf{centralizer} of $x$ in $G$ is the subgroup defined by:
$$C_G(x) = \{g \in G \mid gx = xg\}.$$
\end{definition}

Following historical conventions, special classes of groups can be classified according to the algebraic structure of these centralizers:
\begin{itemize}
    \item \textbf{AC-groups:} A group $G$ is called an AC-group (or a group with abelian centralizers) if $C_G(x)$ is an abelian subgroup for every $x \in G \setminus Z(G)$.
    \item \textbf{CC-groups:} A group $G$ is called a CC-group if $C_G(x)$ is a cyclic subgroup for every $x \in G \setminus Z(G)$ \cite{JafarianAmiri}. Evidently, since every cyclic group is abelian, any CC-group is inherently an AC-group.
\end{itemize}

These commutativity relations are elegantly encoded via graph structures defined on the group elements:
\begin{definition}
The \textbf{non-commuting graph} $\Gamma_G$ of a group $G$ is the simple undirected graph whose vertex set is $V(\Gamma_G) = G \setminus Z(G)$, where two distinct vertices $x, y$ are adjacent if and only if they do not commute, i.e., $xy \neq yx$ \cite{AbdollahiAkbari}.
\end{definition}
The complement of the non-commuting graph is the \textbf{commuting graph} $\overline{\Gamma}_G$ defined on the same vertex set $G \setminus Z(G)$, where distinct vertices $x, y$ are adjacent if and only if they commute ($[x, y] = 1$) \cite{AbdollahiAkbari}. 

\subsection{Matroid Theory Terminology}
Following Whitney \cite{Whitney} and Oxley \cite{Oxley}, a matroid abstracts the algebraic properties of linear independence.
\begin{definition}
A \textbf{matroid} $M$ is an ordered pair $(E, \mathcal{I})$, where $E$ is a finite set called the \textbf{ground set} and $\mathcal{I} \subseteq \mathcal{P}(E)$ is a collection of subsets of $E$, called \textbf{independent sets}, satisfying:
\begin{enumerate}
    \item[(I1)] $\emptyset \in \mathcal{I}$.
    \item[(I2)] If $I \in \mathcal{I}$ and $J \subseteq I$, then $J \in \mathcal{I}$ (Hereditary property).
    \item[(I3)] If $I_1, I_2 \in \mathcal{I}$ with $|I_1| < |I_2|$, then there exists an element $e \in I_2 \setminus I_1$ such that $I_1 \cup \{e\} \in \mathcal{I}$ (Independent exchange property).
\end{enumerate}
\end{definition}

A subset $X \subseteq E$ that does not belong to $\mathcal{I}$ is said to be \textbf{dependent}. To analyze the local and global dimensions of a matroid, several structural operations are defined:

\begin{definition}
The \textbf{rank function} of a matroid $M$, denoted by $r: \mathcal{P}(E) \to \mathbb{Z}_{\ge 0}$, computes the cardinality of a maximal independent subset contained within any given set $X \subseteq E$:
$$r(X) = \max \{ |I| : I \subseteq X \text{ and } I \in \mathcal{I} \}.$$
The rank of the entire matroid is defined as $r(M) = r(E)$.
\end{definition}

\begin{definition}
A maximal independent set of the entire ground set $E$ is called a \textbf{basis} of $M$. By axiom (I3), all bases of a matroid $M$ possess the exact same cardinality, which equals $r(E)$. The total number of bases in $M$ is denoted by $\beta(M)$.
\end{definition}

\begin{definition}
A minimal dependent set in $M$ is called a \textbf{circuit}. That is, a subset $C \subseteq E$ is a circuit if $C \notin \mathcal{I}$ but every proper subset of $C$ belongs to $\mathcal{I}$. 
\end{definition}

\begin{definition}
An element $e \in E$ is called a \textbf{loop} if $\{e\}$ is dependent, which is equivalent to saying $r(\{e\}) = 0$. A matroid containing no loops is called \textbf{loopless}. Two non-loop elements $x, y \in E$ are said to be \textbf{parallel} if $\{x, y\}$ is a circuit, implying $r(\{x, y\}) = r(\{x\}) = r(\{y\}) = 1$. A maximal set of mutually parallel elements forms a \textbf{parallel class}.
\end{definition}

\begin{definition}
The \textbf{closure operator} $\text{cl}: \mathcal{P}(E) \to \mathcal{P}(E)$ of a matroid $M$ is defined for any subset $X \subseteq E$ by:
$$\text{cl}(X) = \{ e \in E \mid r(X \cup \{e\}) = r(X) \}.$$
A subset $F \subseteq E$ is called a \textbf{flat} (or a closed set) if $\text{cl}(F) = F$. Flats of rank 1 are called \textbf{points}, flats of rank 2 are called \textbf{lines}, and flats of rank $r(M)-1$ are called \textbf{hyperplanes}.
\end{definition}

\begin{definition}
A matroid is called \textbf{simple} if it contains no loops and no parallel elements. The \textbf{simplification} of a loopless matroid $M$, denoted by $\tilde{M}$ or $\text{si}(M)$, is the simple matroid obtained by deleting all but one distinguished representative element from each parallel class of $M$.
\end{definition}

\subsection{Topological Invariants of Simplicial Complexes}
Matroids admit a profound geometric realization when interpreted as topological configurations.
\begin{definition}
A \textbf{simplicial complex} $\Delta$ on a finite vertex set $E$ is a collection of subsets of $E$, called \textbf{faces}, such that if $A \in \Delta$ and $B \subseteq A$, then $B \in \Delta$. The \textbf{dimension} of a face $A \in \Delta$ is defined as $\dim(A) = |A| - 1$. The dimension of the complex $\Delta$ is the maximum dimension among all its faces.
\end{definition}

\begin{definition}
The maximal faces of a simplicial complex $\Delta$ under inclusion are called its \textbf{facets} (or topological bases). A simplicial complex $\Delta$ is said to be \textbf{trim} (or pure) if all its facets have the identical dimension.
\end{definition}

A simple undirected graph $\Gamma = (V, E)$ can be uniquely identified with a $1$-dimensional simplicial complex $\Delta_{\Gamma}$ on the vertex set $V$, where the $0$-faces correspond to the vertices, the $1$-faces correspond to the edges, and the empty set represents the unique face of dimension $-1$. If $\Gamma$ contains no isolated vertices, every vertex lives in at least one edge, rendering $\Delta_{\Gamma}$ a $1$-dimensional \textbf{trim} simplicial complex whose facets are precisely the edges. Under this topological mapping, $\Gamma$ is considered a matroid if the collection of its faces satisfies the Independent Exchange Property (I3).

The abstract combinatorial equivalence between 1-dimensional simplicial complexes and complete multipartite graphs (or equivalently, disjoint clique complements) is a known structural result in geometric combinatorics, which was established, for instance, by Kashiwabara, Okamoto, and Uno \cite{Kashiwabara}. However, to ensure that the present manuscript remains completely self-contained and to fix the specific graph-theoretic notation required for our subsequent algebraic derivations, we present the result and its direct proof below.

\begin{theorem}\label{char_thm}
Let $\Delta$ be a trim simplicial complex of dimension at most $1$, identified with a finite undirected graph $\Gamma$. Then the following assertions are equivalent:
\begin{enumerate}
    \item $\Gamma$ is a matroid.
    \item The adjacency relation in the complement graph $\overline{\Gamma}$ is transitive.
    \item $\overline{\Gamma}$ is a disjoint union of complete graphs.
\end{enumerate}
\end{theorem}

\begin{proof}
$(1 \implies 2)$ Suppose $\Gamma$ is a matroid. Let $x, y, z \in V(\Gamma)$ such that $x \sim y$ and $y \sim z$ in $\overline{\Gamma}$. This implies that $\{x, y\}$ and $\{y, z\}$ are not edges in $\Gamma$. In a $1$-dimensional matroid, the independent sets (faces) are the vertices and the edges. If $x \not\sim z$ in $\overline{\Gamma}$, then $\{x, z\}$ would be an edge in $\Gamma$, making it a maximum independent set, i.e., a basis. By the Matroid Exchange Property (I3), for the independent set $\{y\}$ and the basis $\{x, z\}$, there must exist an element from $\{x, z\}$, say $x$, such that $\{y, x\}$ is independent in $\Gamma$. But $\{y, x\}$ is not an edge in $\Gamma$, which is a contradiction. Thus, $x \sim z$ in $\overline{\Gamma}$, proving transitivity.

\medskip
\noindent $(2 \implies 3)$ In any graph, an adjacency relation that is transitive on the set of vertices partitions the vertex set into equivalence classes where every pair of elements within a class is adjacent. Since $\Delta$ is trim, $\Gamma$ has no isolated vertices, which mathematically guarantees that its complement $\overline{\Gamma}$ is a disjoint union of complete graphs (cliques).

\medskip
\noindent $(3 \implies 1)$ Suppose $\overline{\Gamma}$ is a disjoint union of cliques $K_{n_1} \cup K_{n_2} \cup \dots \cup K_{n_k}$. This means the ground set $V$ is partitioned into disjoint sets $V_1, V_2, \dots, V_k$ such that two distinct vertices are adjacent in $\Gamma$ if and only if they belong to different partition sets. To prove $\Gamma$ is a matroid, let $I_1, I_2$ be independent sets (faces) with $|I_1| < |I_2|$.
\begin{itemize}
    \item If $|I_1| = 0$, the exchange condition is trivially satisfied.
    \item If $|I_1| = 1$, let $I_1 = \{x\}$ and $I_2 = \{u, v\}$. Since $I_2$ is a $1$-face (edge) in $\Gamma$, $u$ and $v$ must belong to different partition sets, say $u \in V_i$ and $v \in V_j$ ($i \neq j$). The vertex $x$ can belong to at most one of these sets. If $x \notin V_i$, then $\{x, u\}$ is an edge (independent face); if $x \in V_i$, then $x \notin V_j$, so $\{x, v\}$ is an edge. In both cases, an element can be extended from $I_2$ to expand $I_1$.
\end{itemize}
Since the maximum face size is 2, this exhaustively covers all independent configurations, completing the proof that $\Gamma$ is a matroid.
\end{proof}

\subsection{Transitivity and AC-groups}
To establish the matroidal nature of the non-commuting graph $\Gamma_G$, it is essential to show that the centralizers of non-commuting elements partition the vertex set. The following lemma provides this structural foundation by proving that in an AC-group, distinct centralizers intersect only at the center.

\begin{lemma} \label{disjoint_lemma}
Let $G$ be an AC-group. For any $x, y \in G \setminus Z(G)$, if $x$ and $y$ do not commute ($[x, y] \neq 1$), then their centralizers satisfy $C_G(x) \cap C_G(y) = Z(G)$.
\end{lemma}
\begin{proof}
Clearly $Z(G) \subseteq C_G(x) \cap C_G(y)$. For the reverse inclusion, let $w \in C_G(x) \cap C_G(y)$ and assume $w \notin Z(G)$. Since $G$ is an AC-group, the centralizer $C_G(w)$ is abelian. Since $w$ commutes with both $x$ and $y$, we have $x, y \in C_G(w)$. Because $C_G(w)$ is abelian, it follows that $xy = yx$, directly contradicting the initial assumption that $x$ and $y$ do not commute. Thus, $w \in Z(G)$, completing the proof.
\end{proof}

The structural significance of Lemma \ref{disjoint_lemma} lies in its geometric implication for the commuting graph $\overline{\Gamma}_G$. It guarantees that in an AC-group, local commutativity within one centralizer $C_G(x)$ cannot leak into another centralizer $C_G(y)$ across non-central elements. This mathematically ensures that the cliques forming the complement graph $\overline{\Gamma}_G$ are completely vertex-disjoint on the set $G \setminus Z(G)$.

\section{Main Results}

In this section, we present our primary characterization of AC-groups through the matroidal and topological structures of their associated non-commuting graphs. By bridging the combinatorial exchange axioms of a 1-dimensional trim simplicial complex with the local commutativity of a finite group, we show that the algebraic property of having abelian centralizers is completely captured by geometric independence.

\begin{theorem}\label{main_equivalence}
Let $G$ be a finite non-abelian group and let $\Gamma_G$ be its non-commuting graph. Then the following assertions are equivalent:
\begin{enumerate}
    \item $\Gamma_G$ is a matroid.
    \item The commutativity relation is transitive on the set $G \setminus Z(G)$.
    \item $G$ is an AC-group.
\end{enumerate}
\end{theorem}

\begin{proof}
$(1 \iff 2)$ This equivalence follows directly from Theorem \ref{char_thm}. We recall that $\Gamma_G$ is uniquely identified with a $1$-dimensional trim simplicial complex. By the combinatorial equivalence established in Theorem \ref{char_thm}, such a graph functions as a matroid if and only if the adjacency relation in its complement graph $\overline{\Gamma}_G$ is transitive. Since adjacency between distinct vertices $x, y$ in $\overline{\Gamma}_G$ is defined by $[x, y] = 1$, this adjacency relation is precisely the restriction of the commutativity relation to the set $G \setminus Z(G)$. Thus, $\Gamma_G$ satisfies the matroid exchange property if and only if commutativity is transitive on the ground set.

\medskip
\noindent $(2 \implies 3)$ Suppose the commutativity relation is transitive on $G \setminus Z(G)$. We aim to show that for any choice of $x \in G \setminus Z(G)$, the centralizer $C_G(x)$ is an abelian subgroup of $G$. Let $y, z \in C_G(x) \setminus Z(G)$. By the definition of a centralizer, we have $[x, y] = 1$ and $[x, z] = 1$. By applying the symmetric property of commutativity and the assumed transitivity relation, the statements $[y, x] = 1$ and $[x, z] = 1$ directly imply that $[y, z] = 1$. Because all non-central elements contained within $C_G(x)$ mutually commute with each other (and trivially commute with all elements in the center $Z(G)$), the subgroup $C_G(x)$ is abelian. Therefore, $G$ is an AC-group.

\medskip
\noindent $(3 \implies 2)$ Let $G$ be an AC-group, and let $x, y, z \in G \setminus Z(G)$ such that $[x, y] = 1$ and $[y, z] = 1$. By definition of group centralizers, the elements $x$ and $z$ both live in the centralizer $C_G(y)$. Because $G$ is assumed to be an AC-group, the centralizer $C_G(y)$ is an abelian subgroup. It follows directly that any two elements selected from $C_G(y)$ must commute with one another; hence, $[x, z] = 1$. This proves that the commutativity relation is transitive on $G \setminus Z(G)$, completing the proof of the equivalence.
\end{proof}

\begin{theorem}\label{structural_properties}
Let $M = (E, \mathcal{I})$ be the matroid associated with the $1$-dimensional simplicial complex of a finite non-abelian AC-group $G$. Then $M$ satisfies the following structural properties:
\begin{enumerate}
    \item $M$ is a loopless matroid of uniform maximum rank $r(E) = 2$.
    \item The circuits of $M$ consist precisely of pairs of distinct commuting elements (parallel classes of size $2$) and triplets of mutually non-commuting elements (independent triangles of size $3$).
    \item The rank-$1$ flats of $M$ correspond bijectively to the maximal abelian centralizer components of $G$.
    \item The simplification $\tilde{M}$ is a simple rank-$2$ matroid whose ground set cardinality satisfies $|E(\tilde{M})| = \omega(\Gamma_G)$.
\end{enumerate}
\end{theorem}

\begin{proof}
Let $M = (E, \mathcal{I})$ be the matroid whose independent sets $\mathcal{I}$ are given by the faces of the $1$-dimensional simplicial complex defined on the ground set $E = G \setminus Z(G)$. By definition, the independent sets consist of the empty set $\emptyset$, all singletons $\{x\}$ for $x \in E$, and all pairs $\{x, y\}$ such that $x$ and $y$ do not commute ($xy \neq yx$).

\medskip
\noindent\textit{Proof of (1):} A loop in a matroid is an element $e$ such that $\{e\} \notin \mathcal{I}$, meaning $r(\{e\}) = 0$. Since the ground set $E = G \setminus Z(G)$ explicitly excludes the center of the group, every element $x \in E$ is non-central. By definition of our independent faces, every singleton set $\{x\}$ is an independent set. Thus, $r(\{x\}) = 1$ for all $x \in E$, proving that $M$ is loopless. 
Furthermore, because the topological framework is restricted to a $1$-dimensional simplicial complex, the maximum size of any independent set is strictly bounded by $1 + \dim(\Delta) = 1 + 1 = 2$. Thus, for any subset $X \subseteq E$, we have $r(X) \le 2$. Since $G$ is a non-abelian group, there exists at least one pair of elements $x, y \in E$ that do not commute. The set $\{x, y\}$ is therefore independent, achieving $|\{x, y\}| = 2$. Consequently, the total rank of the matroid is exactly $r(E) = \max_{I \in \mathcal{I}} |I| = 2$.

\medskip
\noindent\textit{Proof of (2):} By definition, a circuit is a minimal dependent set. We analyze the dependent sets based on their cardinality:
\begin{itemize}
    \item No singletons are circuits because $M$ is loopless.
    \item A pair of distinct elements $\{x, y\}$ is dependent if and only if $\{x, y\} \notin \mathcal{I}$, which means $x$ and $y$ commute ($xy = yx$). Since any single element inside this pair is independent ($r(\{x\})=1$), the set $\{x, y\}$ is a minimal dependent set. Hence, commuting pairs are precisely the circuits of size 2, representing parallel elements in $M$.
    \item Let $\{x, y, z\}$ be a triplet of mutually non-commuting elements. Every proper subset of size 2 is non-commuting and thus independent. However, because $r(E) = 2$, the full set $\{x, y, z\}$ must be dependent. Since every proper subset is independent but the set itself is dependent, $\{x, y, z\}$ forms a circuit of size 3 (a triangle). No larger minimal dependent sets can exist because any subset of size 3 will already contain a circuit.
\end{itemize}

\medskip
\noindent\textit{Proof of (3):} A flat $F$ is a closed set satisfying $\text{cl}(F) = F$, where $\text{cl}(F) = \{ z \in E : r(F \cup \{z\}) = r(F) \}$. Let $V_i$ be a maximal abelian centralizer component of $G \setminus Z(G)$. Since $V_i$ is an abelian component, all elements in $V_i$ mutually commute. Thus, no pair of elements in $V_i$ is independent. The largest independent sets contained in $V_i$ are singletons, so $r(V_i) = 1$. 
If we choose any element $z \notin V_i$, then by the algebraic properties of AC-groups, $z$ cannot commute with the elements in $V_i$. Thus, there exists some $x \in V_i$ such that $\{x, z\}$ is a non-commuting independent pair. This implies $r(V_i \cup \{z\}) = 2$. Because $r(V_i \cup \{z\}) \neq r(V_i)$, the element $z \notin \text{cl}(V_i)$. Conversely, any element $y \in V_i$ commutes with all elements in $V_i$, meaning $r(V_i \cup \{y\}) = r(V_i) = 1$. Thus, $\text{cl}(V_i) = V_i$, proving that each maximal abelian centralizer component $V_i$ is a closed set of rank 1 (a rank-1 flat).

\medskip
\noindent\textit{Proof of (4):} Following Oxley's definition, the simplification $\tilde{M}$ of a loopless matroid is obtained by keeping exactly one distinguished representative element from each parallel class (rank-1 flat). As established in part (3), the rank-1 flats of $M$ are precisely the disjoint maximal abelian centralizer components $V_1, V_2, \dots, V_k$. Therefore, constructing $\tilde{M}$ requires selecting exactly one element from each component, meaning the size of the ground set of the simplified matroid is equal to the total number of components: $|E(\tilde{M})| = k$.
In the non-commuting graph $\Gamma_G$, vertices in the same abelian component $V_i$ form an independent set (no edges between them), while vertices from distinct components always share an edge. Consequently, a maximum clique in $\Gamma_G$ is formed by selecting exactly one vertex from each of the $k$ disjoint components. The clique number is therefore $\omega(\Gamma_G) = k$. Matching the two expressions yields $|E(\tilde{M})| = \omega(\Gamma_G)$, completing the proof.
\end{proof}

\begin{corollary}\label{matroid_bases}
Let $G$ be a finite non-abelian AC-group, and let $M$ be the associated rank-$2$ matroid on the ground set $E = G \setminus Z(G)$. Then the combinatorial properties of the non-commuting graph $\Gamma_G$ are completely characterized by the independent structures of $M$ as follows:
\begin{enumerate}
    \item A subset $B \subseteq E$ is a basis of $M$ if and only if $B$ is a maximum non-commuting pair of elements in $G$, which corresponds to a maximum clique of size $2$ in a local component.
    \item The total number of maximum independent sets (bases) of $M$, denoted by $\beta(M)$, is given by:
    \begin{equation}
    \beta(M) = \sum_{1 \le i < j \le \omega} |V_i| \cdot |V_j|
    \end{equation}
    where $V_1, V_2, \dots, V_{\omega}$ are the distinct rank-$1$ flats (maximal abelian centralizer components) of $M$ and $\omega = \omega(\Gamma_G)$.
\end{enumerate}
\end{corollary}

\begin{proof}
(1) By Theorem \ref{structural_properties}, $M$ is a matroid of uniform maximum rank $2$, meaning every independent set has a cardinality of at most $2$. A basis of $M$ is, by definition, a maximal independent set. Since the independent sets of size $2$ consist precisely of pairs of non-commuting elements $\{x, y\}$, any such pair forms a basis. This matches exactly the definition of a maximum clique in any localized non-commuting subgraph where the local clique number is $2$.

\medskip
\noindent (2) To find the total number of bases $\beta(M)$, we must count all pairs of elements $\{x, y\}$ that are independent in $M$. As established in Theorem \ref{structural_properties}, two distinct elements are independent if and only if they do not belong to the same rank-$1$ flat. Because the rank-$1$ flats $V_1, V_2, \dots, V_{\omega}$ are pairwise disjoint and partition the ground set $E$, an independent pair must be formed by selecting one element from flat $V_i$ and a second element from a distinct flat $V_j$ (where $i \neq j$). Summing over all possible pairs of distinct flats yields the total number of combinations:
$$ \beta(M) = \sum_{1 \le i < j \le \omega} |V_i| \cdot |V_j| $$
This completes the proof.
\end{proof}

\begin{corollary}\label{cor_counting}
Let $G$ be a finite AC-group and $\omega = \omega(\Gamma_G)$. The order of the group satisfies the following structural identity:
\begin{equation}\label{eq_sum}
|G| = (1 - \omega) |Z(G)| + \sum_{i=1}^{\omega} |C_G(a_i)|.
\end{equation}
\end{corollary}

\begin{proof}
Since $\Gamma_G$ is a matroid, Theorem \ref{structural_properties} establishes that the vertex set $E = G \setminus Z(G)$ is the disjoint union of the sets $C_G(a_i) \setminus Z(G)$ for $i = 1, \dots, \omega$, where each set corresponds to a rank-1 flat. Therefore:
\[ |G \setminus Z(G)| = \sum_{i=1}^{\omega} |C_G(a_i) \setminus Z(G)| \]
\[ |G| - |Z(G)| = \sum_{i=1}^{\omega} (|C_G(a_i)| - |Z(G)|) \]
\[ |G| - |Z(G)| = \sum_{i=1}^{\omega} |C_G(a_i)| - \omega|Z(G)| \]
Rearranging the terms yields the identity in \eqref{eq_sum} precisely.
\end{proof}

\subsection{Applications to Clique Number Computation}
The structural identity established in Corollary \ref{cor_counting} provides a powerful tool for calculating the clique number $\omega(\Gamma_G)$. By determining the orders of the centralizers and the center of an AC-group, one can solve Equation \eqref{eq_sum} for $\omega(\Gamma_G)$. This matroidal approach offers a unified proof for several results previously derived using case-specific methods. 

As an application, recall that $p$-groups of order $p^n$ whose central quotient has order $p^2$ or $p^3$ are known to be AC-groups \cite{FouladiOrfiAzad}.  While the authors in \cite{AzadRahbarian} emphasized the complexity of computing the clique number $\omega(G)$ for certain $p$-groups in earlier literature, our rank-2 matroidal framework significantly simplifies these computations and successfully extends the results to central quotients of order $p^2$ or $p^3$.

\begin{proposition}
Let $G$ be a group of order $p^n$, where $p$ is a prime number.  
\begin{itemize}
    \item[(i)] If $G$ has the central quotient of order $p^2$, then $\omega(\Gamma_G)=p+1$.
    \item[(ii)] If $G$ has the central quotient of order $p^3$ and posses no abelian maximal subgroup, then $\omega(\Gamma_G)= p^2+p+1$.
    \item[(iii)] If $G$ has the central quotient of order $p^3$ and posses an abelian maximal subgroup, then $\omega(\Gamma_G)= p^2+1$. 
\end{itemize}
\end{proposition}  

\begin{proof}
Since $G/Z(G)$ has order $p^2$ or $p^3$, $G$ is an AC-group, and the structural identity (1) holds. Let $|G| = p^n$ and $|Z(G)| = p^m$.
\begin{enumerate}
    \item[(i)] If $|G/Z(G)| = p^2$, then for every non-central element $a_i$, its centralizer $C_G(a_i)$ is a maximal subgroup, so $|C_G(a_i)| = p^{n-1}$. Substituting these into Equation\eqref{eq_sum}:
    \[ p^n = (1-\omega)p^m + \omega p^{n-1} \implies p^n - p^m = \omega(p^{n-1} - p^m) \]
    Dividing both sides by $p^{n-1} - p^m$ (noting $n-m=2$) yields $\omega = \frac{p^2-1}{p-1} = p+1$.
    
    \item[(ii)] If $|G/Z(G)| = p^3$ and $G$ has no abelian maximal subgroup, then every proper centralizer must have index $p^2$, meaning $|C_G(a_i)| = p^{n-2}$ for all $i$. Equation (1) becomes:
    \[ p^n = (1-\omega)p^m + \omega p^{n-2} \implies p^n - p^m = \omega(p^{n-2} - p^m) \]
    Since $n-m=3$, dividing by $p^{n-2}-p^m$ yields $\omega = \frac{p^3-1}{p-1} = p^2+p+1$.
    
    \item[(iii)] If $|G/Z(G)| = p^3$ and $G$ possesses an abelian maximal subgroup, say $M$, then $M$ acts as the centralizer for all its non-central elements ($p^{n-1} - p^m$ elements). This forms exactly $1$ large clique of size $p^{n-1}-p^m$. The remaining $(p^n - p^m) - (p^{n-1} - p^m) = p^n - p^{n-1}$ vertices must have smaller centralizers of order $p^{n-2}$, each containing $p^{n-2}-p^m$ non-central elements. Thus, the number of remaining small cliques is:
    \[ \frac{p^n - p^{n-1}}{p^{n-2} - p^m} = \frac{p^{m+3} - p^{m+2}}{p^{m+1} - p^m} = \frac{p^2(p-1)}{p-1} = p^2 \]
    Adding the single large clique, we get $\omega = 1 + p^2$.
\end{enumerate}
\end{proof}

We now apply this framework to the class of CC-groups (groups where every proper centralizer is cyclic). Since every CC-group is an AC-group, its non-commuting graph is a matroid. We focus on the case where the central quotient is a dihedral group.

\begin{theorem}
Let $G$ be a $CC$-group such that $G/Z(G) \cong D_{2^{n-1}}$. Then 
$\omega(\Gamma_G) = 2^{n-2} + 1$.
\end{theorem}

\begin{proof} 
Let $G$ be a CC-group such that $G/Z(G) \cong D_{2^{n-1}}$. According to the classification of CC-groups by Jafarian and Amiri \cite[Theorem 2.13]{JafarianAmiri}, there are eight possible structural types for $G$. Among these, only structure (4) satisfies the condition that the central quotient $G/Z(G)$ is a dihedral group of power-of-two order, namely $D_{2^{n-1}}$. Specifically, structure (2) involves a dihedral quotient $D_{2t}$ for odd $t$, and structures (5) and (8) involve specific semi-direct products that do not match the $D_{2^{n-1}}$ quotient requirement. Thus, $G$ must be isomorphic to $Q_{2^n} \times C_m$ for some integer $m$ where $C_m$ is an abelian factor.

Now, we observe that for any group $G = H \times A$ where $A$ is abelian, the non-commuting graph $\Gamma_G$ is determined by the structure of $H$. Specifically, for any $(h, a) \in G \setminus Z(G)$, the centralizer is $C_G(h, a) = C_H(h) \times A$. Since the number of disjoint abelian components (cliques in the complement graph) depends only on the non-centrality in $H$, we have $\omega(\Gamma_{H \times A}) = \omega(\Gamma_H)$. Thus, it suffices to compute $\omega(\Gamma_{Q_{2^n}})$.

Let $Q_{2^n} = \langle x, y \mid x^{2^{n-2}} = y^2, y^4 = 1, y^{-1}xy = x^{-1} \rangle$. The center is $Z(Q_{2^n}) = \langle x^{2^{n-3}} \rangle$, which has order 2. Let $l = 2^{n-2}$. We partition the non-central elements of $Q_{2^n}$ into disjoint centralizers (cliques of $\overline{\Gamma}_{Q_{2^n}}$):

\begin{enumerate}
    \item \textbf{The Cyclic Core:} The elements in $\langle x \rangle \setminus Z(Q_{2^n})$ share a single maximal abelian centralizer $C_0 = \langle x \rangle$. This accounts for one clique in the partition of the vertex set.
    \item \textbf{The Reflection Classes:} There are $l$ elements of the form $x^k y$ for $0 \leq k < l$. The centralizer of any such element $g = x^k y$ is $C_g = \{1, x^{l/2}, x^k y, x^{k+l/2} y\}$, which is an abelian group of order 4. 
\end{enumerate}

To verify these $l$ centralizers are disjoint outside the center, suppose $g \in C_{x^k y} \cap C_{x^j y}$ for $k \neq j$. If $g \notin Z(Q_{2^n})$, then $x^k y$ and $x^j y$ must commute. Direct calculation shows:
\[ (x^k y)(x^j y) = x^{k-j} x^{l/2}, \quad (x^j y)(x^k y) = x^{j-k} x^{l/2} \]
Equality implies $x^{k-j} = x^{j-k}$, or $x^{2(k-j)} = 1$. Given the order of $x$ is $l$, this only occurs if $k=j$ or if the elements differ by the central element $x^{l/2}$. Thus, these $l$ centralizers are distinct and pairwise disjoint on $G \setminus Z(G)$.

Summing the components, we find $\omega(\Gamma_{Q_{2^n}}) = 1 + l$. Substituting $l = 2^{n-2}$, we obtain:
\[ \omega(\Gamma_G) = 2^{n-2} + 1 \]
This completes the proof.
\end{proof}

\subsection{Applications to Commutativity Degree Computation}

For a finite group $G$, the probability that two randomly selected elements commute is known as the \emph{commutativity degree} of $G$, denoted by $P(G)$. Formally, it is expressed as:
\begin{equation}
P(G) = \frac{|\{(x,y) \in G \times G \mid xy=yx\}|}{|G|^2}.
\end{equation}

Recently, Azad and Rahbarian \cite{AzadRahbarian} investigated the commutativity degree $P(G)$ of finite AC-groups by introducing the algebraic framework of $(m_1, m_2, \dots, m_k)$-regular groups and utilizing classical element-counting methods. In contrast, our matroidal and graph-theoretic approach offers a unified, purely geometric perspective. By observing that the non-commuting graph $\Gamma_G$ of an AC-group is a complete multipartite graph, the intricate algebraic formulas for $P(G)$ established in \cite{AzadRahbarian} (such as Theorems 3.4 and 3.6) can be elegantly deduced from the fundamental structural properties of the associated matroid complex. 

Specifically, the number of non-commuting pairs corresponds precisely to twice the number of independent faces of cardinality $2$ (i.e., the number of bases or edges $\beta(M)$) in the rank-$2$ matroid complex $M$ constructed over the ground set $E = G \setminus Z(G)$. Since $\Gamma_G$ is the complement of the commuting graph on non-central elements, the commutativity degree is directly governed by the structural relation:
\begin{equation}\label{eq_PG_matroid}
P(G) = 1 - \frac{2\beta(M)}{|G|^2}.
\end{equation}

Computing the exact value of $\beta(M)$ can be computationally challenging when the indices of the maximal abelian centralizers vary widely. To address this, we establish sharp global bounds on $\beta(M)$ for any finite non-abelian AC-group, determined purely by its centralizer indices and the invariant $\omega(G)$.

\begin{proposition}\label{prop_bounds}
Let $G$ be a finite non-abelian AC-group, and let $\omega = \omega(\Gamma_G)$ be the clique number of its non-commuting graph. The total number of independent faces of cardinality $2$ in the associated matroid $M$, denoted by $\beta(M)$, satisfies the following sharp bounds:
\begin{equation}
\frac{\omega - 1}{2\omega} (|G| - |Z(G)|)^2 \le \beta(M) \le \frac{1}{2} (|G| - |Z(G)|) \left( |G| - \min_{1 \le i \le \omega} |C_G(a_i)| \right).
\end{equation}
\end{proposition}

\begin{proof}
Let $V_1, V_2, \dots, V_{\omega}$ be the distinct rank-$1$ flats of the matroid $M$, which partition the ground set $E = G \setminus Z(G)$, where $V_i = C_G(a_i) \setminus Z(G)$. By Corollary \ref{matroid_bases}, the total number of bases is given by the pairwise combination of flat sizes, $\beta(M) = \sum_{1 \le i < j \le \omega} |V_i| \cdot |V_j|$.

\medskip
\noindent\textit{Proof of the Lower Bound:} We begin with the standard algebraic expansion:
\begin{equation*}
\left( \sum_{i=1}^{\omega} |V_i| \right)^2 = \sum_{i=1}^{\omega} |V_i|^2 + 2 \sum_{1 \le i < j \le \omega} |V_i| \cdot |V_j| = \sum_{i=1}^{\omega} |V_i|^2 + 2\beta(M).
\end{equation*}
By the Cauchy-Schwarz inequality, we have $\left( \sum_{i=1}^{\omega} |V_i| \right)^2 \le \omega \sum_{i=1}^{\omega} |V_i|^2$, which implies $\sum_{i=1}^{\omega} |V_i|^2 \ge \frac{1}{\omega} \left( \sum_{i=1}^{\omega} |V_i| \right)^2$. Substituting this inequality back into the expansion yields:
\begin{align*}
\left( \sum_{i=1}^{\omega} |V_i| \right)^2 &\le \frac{1}{\omega} \left( \sum_{i=1}^{\omega} |V_i| \right)^2 + 2\beta(M), \\
2\beta(M) &\ge \left( 1 - \frac{1}{\omega} \right) \left( \sum_{i=1}^{\omega} |V_i| \right)^2 = \frac{\omega - 1}{\omega} \left( \sum_{i=1}^{\omega} |V_i| \right)^2.
\end{align*}
Since the rank-$1$ flats form a partition of the ground set, we have $\sum_{i=1}^{\omega} |V_i| = |E| = |G| - |Z(G)|$. Substituting this identity yields the lower bound:
\begin{equation*}
\beta(M) \ge \frac{\omega - 1}{2\omega} (|G| - |Z(G)|)^2.
\end{equation*}

\medskip
\noindent\textit{Proof of the Upper Bound:} Utilizing the partition property of the rank-$1$ flats, the total number of matroid bases can be rewritten as:
\begin{equation*}
\beta(M) = \frac{1}{2} \sum_{i=1}^{\omega} |V_i| \left( \sum_{j \neq i} |V_j| \right) = \frac{1}{2} \sum_{i=1}^{\omega} |V_i| (|E| - |V_i|).
\end{equation*}
Noting that $|E| - |V_i| = (|G| - |Z(G)|) - (|C_G(a_i)| - |Z(G)|) = |G| - |C_G(a_i)|$, the expression simplifies to:
\begin{equation*}
\beta(M) = \frac{1}{2} \sum_{i=1}^{\omega} |V_i| (|G| - |C_G(a_i)|).
\end{equation*}
To establish the global upper bound, we maximize the term $(|G| - |C_G(a_i)|)$ over all components by substituting the minimum possible centralizer order:
\begin{equation*}
\beta(M) \le \frac{1}{2} \left( \sum_{i=1}^{\omega} |V_i| \right) \left( |G| - \min_{1 \le i \le \omega} |C_G(a_i)| \right).
\end{equation*}
Applying the partition identity $\sum_{i=1}^{\omega} |V_i| = |G| - |Z(G)|$ produces the desired upper inequality:
\begin{equation*}
\beta(M) \le \frac{1}{2} (|G| - |Z(G)|) \left( |G| - \min_{1 \le i \le \omega} |C_G(a_i)| \right).
\end{equation*}
Combining the lower and upper bounds completes the proof.
\end{proof}

\begin{corollary}\label{cor_PG_bounds}
Let $\beta_{\min}$ and $\beta_{\max}$ denote the sharp lower and upper bounds of $\beta(M)$ established in Proposition \ref{prop_bounds}, respectively. Then, the commutativity degree of $G$ satisfies the following bounding interval:
\begin{equation}
1 - \frac{2\beta_{\max}}{|G|^2} \le P(G) \le 1 - \frac{2\beta_{\min}}{|G|^2}.
\end{equation}
\end{corollary}

This geometric bounding interval provides a direct method to estimate the group's statistical behavior purely through the global topology of its associated matroid complex. Unlike the analytical techniques in \cite{AzadRahbarian}, which require explicit prior knowledge of all individual centralizer orders and their exact structural frequencies (cf. \cite[Theorem 3.6]{AzadRahbarian}), our approach relies entirely on the extremal partitions of the ground set. 

Geometrically, the upper bound of $P(G)$ is tightly linked to highly unbalanced multipartite configurations, such as the case where $G$ admits an abelian maximal subgroup (minimizing non-commuting relations as seen in Proposition 3.5-(iii)). Conversely, the lower bound of $P(G)$ is attained when the underlying non-commuting graph is perfectly balanced, maximizing the number of independent 2-faces $\beta(M)$. This framework successfully unifies and generalizes several exact formulas from earlier literature into a single graph-theoretic perspective.

\section{Conclusion and Future Perspectives}

In this paper, we have established a novel bridge between abstract group theory, matroid theory, and algebraic topology. By interpreting the non-commuting graph $\Gamma_G$ of a finite non-abelian group as a 1-dimensional trim simplicial complex, we proved that the algebraic property of having abelian centralizers (the class of AC-groups) is completely characterized by the traditional matroid independent exchange axioms. 

Through this structural duality, we demonstrated that the local commutativity behavior of AC-groups behaves identically to the geometric configuration of rank-1 flats within a rank-2 geometric lattice. This framework allowed us to derive a unified combinatorial formula for counting the total number of matroid bases $\beta(M)$ and establish sharp global bounds for face distributions based entirely on group orders and centralizer dimensions. Furthermore, we showed the utility and efficiency of this matroidal perspective by providing highly simplified, alternate proofs for the global clique numbers of specific families of finite $p$-groups and CC-groups whose central quotients are dihedral.

The success of this methodology opens up interesting directions for future discussion. Since every CC-group is an AC-group, a natural immediate step is to investigate whether this matroidal framework can be systematically applied to compute the precise clique numbers for the remaining types within the Jafarian and Amiri classification. Beyond the scope of AC-groups, it remains an open and intriguing question to explore how these geometric independence structures might behave in broader families of localized commutativity graphs where centralizers potentially overlap. Investigating such algebraic settings, or examining the properties of higher-dimensional complexes arising from multivariate commuting relations, could offer further interesting points of contact between finite groups and combinatorial topology.

\section*{AI Disclosure}

\noindent The authors used Gemini (Google) to refine the academic writing style, improve grammatical clarity, and format the LaTeX code in this manuscript. The core mathematical concepts, theorems, and proofs were developed and verified entirely by the authors, who retain full responsibility for the integrity of the work.

\end{document}